# Estimation of the asymptotic behavior of summation functions

Victor Volfson

ABSTRACT  The paper considers asymptotics of summation functions of additive and multiplicative arithmetic functions. We also study asymptotics of summation functions of natural and prime arguments. Several assertions on this subject are proved and examples are considered.

# 1. INTRODUCTION

An arithmetic function is a function defined on the set of natural numbers and taking values on the set of complex numbers. The name arithmetic function is due to the fact that this function expresses some arithmetic property of the natural series.

A summation function is a function of the form:

$$S(n) = \sum_{m \leq n} f(m), \qquad (1.1)$$

where $f(m), m = 1,...,n$ is an arithmetic function.

The Mertens function $M(n) = \sum_{m \leq n} \mu(m)$ is an example of a summation function, where $\mu(m)$ is the arithmetic Möbius function. The value of Möbius function $\mu(m) = 1$ if the natural number m has an even number of prime divisors of the first degree, $\mu(m) = -1$ if the natural number m has an odd number of prime divisors of the first degree, and $\mu(m) = 0$ if the natural number m has prime divisors not only of the first degree.

The average value of the arithmetic function is determined by the formula:

$$E[f,n] = \frac{\sum_{m \leq n} f(m)}{n}. \qquad (1.2)$$

Based on (1.1), (1.2), the asymptotic of the summation function at $n \to \infty$ is determined by the formula:

$$S(n) = E[f,n]n. \qquad (1.3)$$

Thus, based on (1.3), obtaining the asymptotic of the summation function is reduced to obtaining the asymptotic of the mean value of the arithmetic function $f(m), m = 1,...,n$.

The density of the summation function is determined by the formula:

$$d(S) = \frac{S(n)}{n}. \qquad (1.4)$$

The asymptotic density of the summation function is the limit:

$$d^*(S) = \lim_{n \to \infty} \frac{S(n)}{n}. \tag{1.5}$$

Proving the existence of this limit (the asymptotic density of the summation function) for various arithmetic functions is an important problem in finding asymptotics of summation functions.

Let there are two coprime natural numbers $m_1, m_2$. Then an arithmetic function $f$ is called additive if:

$$f(m_1 m_2) = f(m_1) + f(m_2). \tag{1.6}$$

An arithmetic function $f$ is called multiplicative if:

$$f(m_1 m_2) = f(m_1) f(m_2). \tag{1.7}$$

The arithmetic function $f(m), m = 1,...,n$ (in the summation function) can be multiplicative, additive, or some other.

An arithmetic function $f(m), m = 1,...,n$ has a normal order about its mean [1] if:

$$f(n) = \mathrm{E}[f, n](1 + o(1)). \tag{1.8}$$

Obtaining the asymptotic of the summation function has been and is a topical problem in probabilistic number theory. However, summation functions with a multiplicative arithmetic function are studied in most papers [2], [3], [4], [5]. This work should partially compensate for this.

## 2. ASYMPTOTIC BEHAVIOR OF SUMMATION FUNCTIONS WITH ADDITIVE ARITHMETIC FUNCTION

Let a natural number $m$ have a canonical decomposition $m = p_1^{a_1}...p_t^{a_t}$, where $p_i$ is a prime number and $\alpha_i$ is a natural number.

Then, the property for the additive arithmetic function is fulfilled:

$$f(m) = f(p_1^{a_1}...p_t^{a_t}) = f(p_1^{a_1}) + ... + f(p_t^{a_t}) = \sum_{p^\alpha | m} f(p^\alpha). \tag{2.1}$$

A strongly additive arithmetic function is a function for which $f^*(p^\alpha) = f(p)$. Therefore, based on (2.1), the following property holds for a strongly additive arithmetic function:

$$f^*(m) = f^*(p_1^{a_1} \ldots p_t^{a_t}) = \sum_{p|m} f(p). \qquad (2.2)$$

First, we consider the asymptotic of the mean value for a strongly additive arithmetic function $f^*(m) = \sum_{p|m} f(p)$.

Assertion 1

The asymptotic of the mean value of a strongly additive arithmetic function $f^*(m)$ on the interval $[1, n]$ at value $n \to \infty$ is equal to:

$$A_n = \sum_{p \leq n} \frac{f(p)}{p}. \qquad (2.3)$$

Proof

We define a random variable $f^{(p)}(m) = f(p)$ for each prime $p(p \leq n)$ $f^{(p)}(m) = f(p)$ if $p|m$ and $f^{(p)}(m) = 0$ otherwise.

Then $f^{(p)}(m) = f(p)$ with probability $\frac{1}{n}[\frac{n}{p}]$ and $f^{(p)}(m) = 0$ with probability $1 - \frac{1}{n}[\frac{n}{p}]$.

Therefore, the average value $f^{(p)}(m)$ over the interval $[1, n]$ is:

$$E[f^{(p)}, n] = \frac{f(p)}{n}[\frac{n}{p}].$$

It is performed for a strongly additive arithmetic function $f^*(m)$:

$$f^*(m) = \sum_{p \leq n} f^{(p)}(m). \qquad (2.4)$$

Based on (2.4), the average value $f^*(m)$ on the interval $[1, n]$ is:

$$E[f^*, n] = \sum_{p \leq n} \frac{f(p)}{n}[\frac{n}{p}].$$

Therefore, the desired asymptotic of the mean value of a strongly additive arithmetic function $f^*(m)$ on the interval $[1,n]$ with the value $n \to \infty$ is:

$$A_n = \sum_{p \leq n} \frac{f(p)}{p}.$$

The class S of additive arithmetic functions was introduced in [6] for which the additive arithmetic function $f(m), m=1,...,n$ and the corresponding strongly additive arithmetic function $f^*(m) = \sum_{p|m} f(p)$ have the same asymptotic behavior of the probability characteristics for the value $n \to \infty$.

The assertion is proved in the same work that if the condition is satisfied for the additive arithmetic function $f(m), m=1,...,n$ with the value $n \to \infty$:

$$f(n) = O(\ln(n)), \qquad (2.5)$$

then the additive arithmetic function $f(m)$ and the corresponding strongly additive arithmetic function $f^*(m) = \sum_{p|m} f(p)$ belong to the class $S$.

Thus, when condition (2.5) is satisfied, formula (2.3) can be used to find the asymptotic of the mean value of an additive arithmetic function at $n \to \infty$.

Let's look at an example of finding the asymptotic mean of an additive arithmetic function $f(m) = \ln \varphi(m)$. Let us first show that the additive arithmetic function $f(m) = \ln \varphi(m)$ satisfies condition (2.5).

It is known that:

$$\varphi(m) = m \prod_{p|m} (1-1/p) \leq m. \qquad (2.6)$$

Based on (2.6):

$$\ln \varphi(n) = O(\ln(n)). \qquad (2.7)$$

It follows from (2.7) that the additive arithmetic function $f(m) = \ln \varphi(m)$ satisfies the conditions of the assertion [6] and belongs to the class $S$, i.e. the asymptotics of the probabilistic characteristics of this function coincide with the asymptotics of the probabilistic characteristics of a strongly additive arithmetic function $f^*(m) = \sum_{p|m} f(p) = \sum_{p|m} \ln \varphi(p)$.

Having in mind $\varphi(m) = m \prod_{p|m}(1-1/p)$ we get:

$$f(p) = \ln \varphi(m) = \ln(p(1-1/p)) = \ln p + \ln(1-1/p). \tag{2.8}$$

Based on (2.3) and (2.8), the asymptotic of the mean value of a strongly additive arithmetic function $f^*(m), m = 1,...,n$ at value $n \to \infty$ is equal to:

$$\sum_{p \leq n} \frac{f(p)}{p} = \sum_{p \leq n} \frac{\ln p}{p} + \sum_{p \leq n} \frac{\ln(1-1/p)}{p}. \tag{2.9}$$

Considering that based on [7] the asymptotic $\sum_{p \leq n} \frac{\ln p}{p} = \ln n + O(1)$ is fulfilled,

and the series $\sum_{p} \frac{\ln(1-1/p)}{p}$ - converges, then based on (2.9) we obtain the asymptotic of the mean value of a strongly additive arithmetic function $f^*(m), m = 1,...,n$ for the value $n \to \infty$:

$$A_n = \ln n + O(1). \tag{2.10}$$

Based on the assertion [6], the asymptotic of the mean value of an additive arithmetic function $f(m) = \ln \varphi(m)$ is also determined by formula (2.10).

Taking into account (2.3) and (2.10), the asymptotic of the summation function at value $n \to \infty$ is equal to:

$$\sum_{m \leq n} \ln(\varphi(m)) = n \ln n + O(n). \tag{2.11}$$

Let's consider one more example. It is required to find the asymptotic of the summation functions $\sum_{m \leq n} \omega(m)$ and $\sum_{m \leq n} \Omega(m)$.

First, based on (2.4), we determine the asymptotic of the mean value of a strongly additive arithmetic function $\omega(m), m = 1,...,n$ for the value $n \to \infty$:

$$\sum_{p \leq n} \frac{\omega(p)}{p} = \sum_{p \leq n} \frac{1}{p} = \ln \ln n + O(1). \tag{2.12}$$

It was shown in [6] that the conditions of the assertion are satisfied for additive arithmetic functions $\omega(m), \Omega(m)$. Therefore, the asymptotic of the mean value of the additive arithmetic function $\Omega(m), m = 1,...,n$, with the value $n \to \infty$, is also determined by the formula (2.12).

Based on (2.3) and (2.12), asymptotics of the indicated summation functions at value $n \to \infty$ is equal to:

$$\sum_{m \leq n} \Omega(m) = \sum_{m \leq n} \omega(m) = n \ln \ln n + O(n). \tag{2.13}$$

## 3. ESTIMATION OF THE ASYMPTOTIC BEHAVIOR OF SUMMATION FUNCTIONS WITH MULTIPLICATIVE ARITHMETIC FUNCTIONS

Most theorems about the mean value of a multiplicative arithmetic function $g(m), m = 1,...,n$ consider the case when $|g(m)| \leq 1$.

For example, Wirzing's theorem [3]: Let there is a real multiplicative function $g$ with value in $[-1,1]$. Then it executes:

$$lim_{n \to \infty} \frac{1}{n} \sum_{m \leq n} g(m) = \prod_{p} (1 - 1/p) \sum_{\nu \geq 0} \frac{g(p^\nu)}{p^\nu}$$

where the infinite product should be taken equal to zero when it diverges.

Delange's theorem [4]: Let g be a complex-valued multiplicative function with values in the unit circle (or real with module not exceeding 1). Then, provided that $\sum_{p} \frac{1 - \operatorname{Re} g(p)}{p} < \infty$, we get:

$$\sum_{m \leq n} g(m) = n \prod_{p \leq n} (1 - 1/p) \sum_{\nu \geq 0} \frac{g(p^\nu)}{p^\nu} + o(n) \quad (n \to \infty).$$

Kubilius' theorem [5]. Let g be a complex-valued multiplicative function, $|g(m)| \leq 1$. Assume that there is a quantity $\kappa$ independent of $p$ and a constant $c_s$, satisfying the inequality:

$$\sum_p |g(p)-\kappa|\frac{\ln p}{p} < c_s$$

Then at $n \geq 20$:

$$\sum_{m \leq n} g(m) = \frac{n(\ln n)^{\kappa-1}}{\Gamma(\kappa)} \prod_p (1-\frac{1}{p})^\kappa (1+\sum_{\alpha=1}^{\infty} \frac{g(p^\alpha)}{p^\alpha}) + Bn\sqrt{\frac{\ln \ln n}{\ln n}}.$$

Here, $\Gamma(\kappa)$ is gamma function. The infinite product converges absolutely. The constant $B$ depends only on $c_s$. It's obvious that $|\kappa| \leq 1$. We consider that the value $1/\Gamma(\kappa) = 0$ in cases $\kappa = 0$ and $\kappa = -1$.

Considering the above, it is interesting to consider the asymptotic of the mean value for the multiplicative function in other cases.

Assertion 2

Let $g(m), m = 1,...,n$ is a real multiplicative arithmetic function that has normal order with respect to its mean value and $g(m) > 0$. Then the asymptotic of the mean value $g(m), m = 1,...,n$ almost everywhere at value $n \to \infty$ is determined by the formula:

$$\frac{1}{n}\sum_{m \leq n} g(m) = e^{E[f,n]+O(b(n)\sqrt{D[f,n]})}(1+o(1)),$$

where $E[f,n]$, $D[f,n]$ are respectively the mean value and variance of the corresponding additive arithmetic function $f(m), m = 1,...,n$ at the value $n \to \infty$, and $b(n)$ is a slowly growing function.

Proof

Let a natural number $m$ have a canonical decomposition $m = p_1^{a_1}...p_t^{a_t}$, where $p_i$ is a prime number and $\alpha_i$ is a natural number.

A multiplicative arithmetic function $g(m), m = 1,...,n$ is converted to an additive one using the logarithm operation:

$$f(m) = \ln(g(m)) = \ln(g(p_1^{a_1}...p_t^{a_t})) = \ln(\prod_{p^\alpha | m} g(p)) = \sum_{p^\alpha | m} \ln(g(p)). \quad (3.1)$$

It was proved in [8] that any arithmetic function $f(m), m = 1,...,n$ almost everywhere has the following asymptotic for the value $n \to \infty$:

$$f(n) = E[f,n] + O(b(n)\sqrt{D[f,n]}), \qquad (3.2)$$

where $E[f,n]$, $D[f,n]$ are respectively the mean value and variance $f(m), m = 1,...,n$ at the value $n \to \infty$, and $b(n)$ is a slowly growing function.

Having in mind that $f(m), m = 1,...,n$ is an additive arithmetic function, then, based on (3.1) and (3.2), the asymptotic behavior $g(m), m = 1,...,n$ almost everywhere at $n \to \infty$ is equal to:

$$g(n) = e^{E[f,n] + O(b(n)\sqrt{D[f,n]})}. \qquad (3.3)$$

Since, by the condition, the multiplicative arithmetic function $g(m), m = 1,...,n$ has a normal order with respect to its mean value, then, based on (1.8) and (3.3), we obtain the desired asymptotic almost everywhere for the mean value $g(m), m = 1,...,n$ at the value $n \to \infty$:

$$\frac{1}{n}\sum_{m \leq n} g(m) = e^{E[f,n] + O(b(n)\sqrt{D[f,n]})}(1 + o(1)). \qquad (3.4)$$

Let's look at an example for assertion 2.

Let it is required to determine the asymptotic almost everywhere for the strongly multiplicative arithmetic function:

$$g^*(m) = e^{\omega(m) + \sum_{p|m} \ln(1 - 1/p)}, m = 1,...,n \text{ at } n \to \infty. \qquad (3.5)$$

First, we note that the strongly multiplicative arithmetic function $g^*(m)$ is not bounded as $m$ grows.

The asymptotic almost everywhere for the strongly multiplicative arithmetic function (3.5) was found in [8] for the value $n \to \infty$:

$$g^*(n) = \ln(n) e^{O((\ln\ln n)^{1/2 + \xi})}. \qquad (3.6)$$

It is also shown there that the strongly multiplicative arithmetic function (3.5) has a normal order with respect to its mean value. Consequently, the conditions of Asertion 2 are satisfied for the strongly multiplicative arithmetic function (3.5).

Therefore, based on (3.4) and (3.6), the asymptotic almost everywhere for this strongly multiplicative arithmetic function at value $n \to \infty$ is equal to:

$$\frac{1}{n}\sum_{m\leq n} g^*(m) = \ln(n)e^{O((\ln\ln n)^{1/2+\xi})}(1+o(1)). \tag{3.7}$$

Based on (1.3) and (3.7), the asymptotic almost everywhere for the summation function of a given multiplicative arithmetic function for $n \to \infty$ will be equal to:

$$\sum_{m\leq n} g^*(m) = n\ln(n)e^{O((\ln\ln n)^{1/2+\xi})}(1+o(1)). \tag{3.8}$$

## 4. ESTIMATION OF THE ASYMPTOTIC BEHAVIOR OF SUMMATIONS FUNCTIONS OF THE NATURAL ARGUMENT

Assertion 3

Let there is a limit $d^*(S) = \lim_{n\to\infty} \frac{S(n)}{n}$, then the asymptotic of the summation function $S(n) = \sum_{m\leq n} f(m)$ for $n \to \infty$ has the form:

$$S(n) = d^*(S)n + o(n), \tag{4.1}$$

Proof

If $d^*(S) = \lim_{n\to\infty} \frac{S(n)}{n} = 0$, then:

$$S(n) = o(n), \tag{4.2}$$

which corresponds to (4.1).

If $d^*(S) \neq 0$, then we get:

$$\lim_{n\to\infty}\frac{S(n) - d^*(S)n}{n} = \left(\lim_{n\to\infty}\frac{S(n)}{n}\right) - d^*(S) = 0. \tag{4.3}$$

Based on (4.3) we get:

$$S(n) - d^*(S)n = o(n).$$

This means that the desired asymptotic also holds in this case:

$$S(n) = d^*(S)n + o(n).$$

An example of the fulfillment of assertion 3 is the asymptotic of the summation function of the number of natural numbers belonging to any subset of natural numbers $A$:

$$\pi(A, n) = |\{m \leq n : m \in A\}| \leq n.$$

So in this case:

$$d^*(\pi) = \lim_{n \to \infty} \frac{\pi(A, n)}{n} \leq 1.$$

An example of this case is the asymptotic of the number of natural numbers belonging to the subset of square-free numbers - Q:

$$\pi(Q, n) = d^*(\pi)n + o(n), \tag{4.4}$$

where $d^*(\pi) = \dfrac{6}{\pi^2} < 1$.

Another example is the asymptotic of the number of natural numbers belonging to a subset of primes - P:

$$\pi(P, n) = o(n), \tag{4.5}$$

where $d^*(\pi) = 0$.

Assertion 4

If $f(m)$ is a multiplicative function and $|f(m)| \leq 1$, then the next asymptotic is true for the summation function:

$$S(n) = \sum_{m \leq n} f(m) = n \prod_{p} (1 - 1/p) \sum_{\nu \geq 0} \frac{f(p^\nu)}{p^\nu} + o(n), \tag{4.6}$$

where the infinite product is considered to be 0 if the product diverges.

The proof follows from Wirzing's theorem [3] and Assertion 3.

An example of the fulfillment of Assertion 4 is the asymptotic of the Mertens summation function $M(n) = \sum_{m \leq n} \mu(m)$, which we talked about in the introduction. The Möbius function is multiplicative and $|\mu(m)| \leq 1$.

In this case, the infinite product diverges, so:

$$d^*(M) = 0. \tag{4.7}$$

and based on (4.6) and (4.7), the Mertens function satisfies the asymptotic:

$$M(n) = o(n). \tag{4.8}$$

Let's look at the summation function $S(n) = \sum_{m \leq n} \frac{1}{m}$. This arithmetic function is also a multiplicative one, and $|f(m)| \leq 1$ for it, i.e. it satisfies the conditions of Assertion 4 and, based on this assertion, the asymptotic of its summation function is also equal to:

$$\sum_{m \leq n} \frac{1}{m} = o(n). \tag{4.9}$$

However, the arithmetic function $f(m) = 1/m$ is elementary and a more precise asymptotic estimate is true for it.

In the general case, if an arithmetic function $f$ is elementary and sufficiently smooth on the interval $[1, n]$, then it is better to use the Euler-Macleron formula [9] to determine the asymptotic of its summation function.

In particular, if an arithmetic function $f$ is strictly decreasing on the interval $[1, n]$ and $\lim_{n \to \infty} f(n) = 0$, then the asymptotic behavior of its summation function is determined by the formula:

$$S(n) = \sum_{m \leq n} f(m) = \int_{t=1}^{n} f(t)dt + O(1). \tag{4.10}$$

As an example of this case, we will consider the summation function indicated above $S(n) = \sum_{m \leq n} \frac{1}{m}$. Based on (4.10), this summation function has the asymptotic:

$$\sum_{m \leq n} \frac{1}{m} = \ln(n) + O(1). \tag{4.11}$$

If we compare the asymptotic estimates (4.9) and (4.11) for this summation function, then estimate (4.11) is more accurate.

If the arithmetic function $f$ is non-decreasing on the interval $[1, n]$, then the asymptotic behavior of the summation function is determined by the formula:

$$S(n) = \sum_{m \leq n} f(m) = \int_{t=1}^{n} f(t)dt + O(f(n)). \tag{4.12}$$

As an example of using (4.12), we consider the definition of the asymptotic of the following summation function:

$$S(n) = \sum_{m=1}^{n} m^k = \int_{t=1}^{n} t^k dt + O(n^k) = \frac{n^{k+1}}{k+1} + O(n^k), \tag{4.13}$$

where $k > 0$ is a constant.

The asymptotic density of the summation function in this case is equal to:

$$d^*(S) = \lim_{n \to \infty} \frac{S(n)}{n} = \lim_{n \to \infty} \frac{\frac{n^{k+1}}{k+1} + O(n^k)}{n} = \infty,$$

so assertion 3 cannot be used.

5. ESTIMATION OF THE ASYMPTOTIC BEHAVIOR OF SUMMATION FUNCTIONS OF THE PRIME ARGUMENT

Assertion 5

If $f'(t)$ exists and is continuous, then:

$$\sum_{p \leq n} f(p) = \frac{nf(n)}{\ln(n)} + O(\frac{n|f(n)|}{\ln^2(n)}) - \int_2^n \frac{tf'(t)dt}{\ln(t)} + O(\int_2^n \frac{t|f'(t)|dt}{\ln^2(t)}) \tag{5.1}$$

Proof

Let $a_k = 1$ if $k$ is a prime number and $a_k = 0$ otherwise.

Let us denote $A(n) = \sum_{k=1}^{n} a_k = \pi(n)$, where $\pi(n)$ is the number of prime numbers not exceeding the value $n$.

If $f'(t)$ exists and is continuous, then using the Abel formula to the sum $\sum_{k=1}^{n} a_k f(k)$, we get:

$$\sum_{p \leq n} f(p) = \sum_{k=1}^{n} a_k f(k) = A(n) f(n) - \int_{1}^{n} A(t) f'(t) dt . \qquad (5.2)$$

Based on the asymptotic law of prime numbers:

$$A(n) = \pi(n) = \frac{n}{\ln(n)} + O(n / \ln^2(n)) . \qquad (5.3)$$

Now we substitute (5.3) into (5.2) and get:

$$\sum_{p \leq n} f(p) = \frac{n f(n)}{\ln(n)} + O(\frac{n |f(n)|}{\ln^2(n)}) - \int_{2}^{n} \frac{t f'(t) dt}{\ln(t)} + O(\int_{2}^{n} \frac{t |f'(t)| dt}{\ln^2(t)}) ,$$

which corresponds to assertion 5.

Let us give an example of using formula (5.1) to obtain the asymptotic of the Chebyshev function:

$$\upsilon(n) = \sum_{p \leq n} \ln(p) = n + O(n / \ln(n)) . \qquad (5.4)$$

Based on (5.4), the asymptotic density of the Chebyshev summation function is:

$$d^*(\upsilon) = \lim_{n \to \infty} \frac{n + O(n / \ln(n))}{n} = 1 . \qquad (5.5)$$

Since the conditions of Assertion 3 are satisfied for the Chebyshev summation function, then based on (4.1) and (5.5) we obtain:

$$\upsilon(n) = n + o(n). \qquad (5.5)$$

The asymptotic estimate (5.4) is more accurate than the estimate (5.5).

Therefore, if an arithmetic function $f$ is a function of the prime argument and satisfies the conditions of Assertion 5, then it is better to use formula (5.1) to determine the asymptotic of the summation function.

## 6. CONCLUSION AND SUGGESTIONS FOR FURTHER WORK

The next article will continue to study the asymptotic behavior of some arithmetic functions.

## 7. ACKNOWLEDGEMENTS

Thanks to everyone who has contributed to the discussion of this paper. I am grateful to everyone who expressed their suggestions and comments in the course of this work.